\documentclass{article}
\usepackage{amsmath,amssymb,graphicx}
\newtheorem{df}{Definition}[section]
\newtheorem{thm}[df]{Theorem}
\newtheorem{lem}[df]{Lemma}

\title{Angular multiselectivity with spherical wavelets}
\author{Ilona Iglewska-Nowak\footnote{West Pomeranian University of Technology in Szczecin, School of Mathematics, al. Pias\-t\'ow 17, 70--310 Szczecin, Poland}}

\begin{document}

\maketitle

\bibliographystyle{amsplain}

\begin{abstract}We construct spherical wavelets based on approximate identities that are directional, i.e. not rotation-invariant, and have an adaptive angular selectivity. The problem of how to find a proper representation of distinct kinds of details of real images, ranging from highly directional to fully isotropic ones, was quite intensively studied for the case of signals over the Euclidean space. However, the present paper is the first attempt to deal with this task in the case of spherical signals. A multiselectivity scheme, similar to that proposed for $\mathbb R^2$-functions, is presented.
\end{abstract}

\begin{bfseries}Key words and phrases:\end{bfseries} spherical wavelets, Poisson kernel, angular selectivity, directional wavelet frames, angular multiresolution \\
\begin{bfseries}2010 Mathematics Subject Classification Number:\end{bfseries} 42C40, 42C15, 94A12, 65T60

\section{Introduction}

Recently, directional wavelets based on approximate identities have been introduced \cite{sB09,EBCK09,HH09,IIN15CWT,IIN16DW}. Unfortunately, wavelet families constructed so far have a very low angular resolution, e.g., the first directional derivative of the Poisson kernel is a wavelet with separable spherical variables $(\vartheta,\varphi)$ of the form
$$
\Psi_\rho^{[1]}(x)=\psi_\rho^{[1]}(\vartheta)\cdot\cos\varphi,
$$
where
$$
x=(\cos\vartheta,\sin\vartheta\cos\varphi,\sin\vartheta\sin\varphi)\in\mathcal S,\qquad\vartheta\in[0,\pi],\,\varphi\in[0,2\pi).
$$
(Unless it leads to misunderstandings, we identify $x\in\mathcal S$ with its spherical coordinates $(\vartheta,\varphi)$.)
\begin{figure}\label{fig:two_wavelets}
\centering\vspace{-8em}\hspace{2em}
\includegraphics[angle=270,scale=0.4]{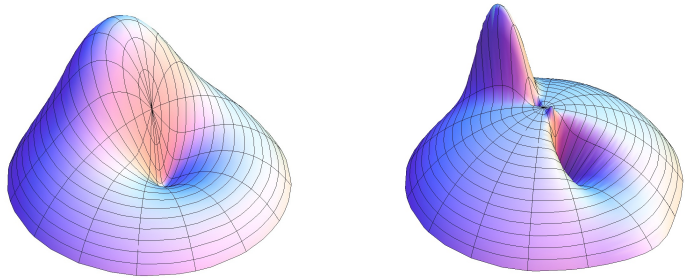}\vspace{-7em}
\caption{The first directional derivative of the Poisson kernel on $\mathcal S^2$ vs. a wavelet with a bigger angular selectivity -- behavior around the north pole}
\end{figure}It is depicted on the left-hand side of Figure~\ref{fig:two_wavelets} and it can be seen that its ability to detect features with close orientations is not satisfactory. In~\cite{HH09} second order wavelets are used to analyse the land cover data of the world, where the analyzed objects are in principle curves. The wavelets are of the form
$$
\Psi_\rho^{[2]}(x)=\psi_{\rho,1}^{[2]}(\vartheta)+\psi_{\rho,2}^{[2]}(\vartheta)\cos^2\varphi,
$$
see \cite[Appendix~A]{IIN16DW}, i.e. again functions with a low angular resolution. The aim of the present paper is to construct wavelets with a variable angular selectivity such that their ability to detect directional features can be adapted to the analyzed signal. An example of such a wavelet is presented on the right-hand-side of Figure~\ref{fig:two_wavelets}. It is a function with separable variables and Figure~\ref{fig:f_tau} shows the dependence of the function of the longitudinal variable on a parameter that can be chosen almost freely.

The investigation is inspired by the results obtained by Antoine and Jacques in \cite{AJ03book,AJ03,JA07} for the case of wavelets over~$\mathbb R^2$. The angular selectivity of spherical wavelets was analysed by~\cite{AMVA08}, but in this case wavelets based on the group theory were involved (the difference between the two constructions is discussed in~\cite{IIN15CWT}). The angular resolving power was defined in \cite[Subsection~9.2.4.2]{AMVA08}. Nonetheless, we do not intend to adapt it to our purposes, since we find the parameter~$\tau$ introduced in Section~\ref{sec:construction} to be a better indicator of the angular selectivity in the case of wavelets constructed in the present paper. Moreover, one that can be a priori chosen.

The paper is organized as follows. After an introduction of necessary notions in Section~\ref{sec:sphere}, we present a construction of two wavelets with a steerable angular selectivity in Section~\ref{sec:construction}. An angular multiselectivity analysis is then discussed in Section~\ref{sec:multiselectivity}.

\section{Preliminaries}\label{sec:sphere}

By~$\mathcal S$ we denote the unit two-dimensional sphere in~$\mathbb R^3$, \mbox{$\mathcal S=\{\xi\in\mathbb R^3:\,|\xi|^2=1\}$}, with the
measure~$d\sigma$ invariant under the rotation group and such that \mbox{$\int_\mathcal S d\sigma=4\pi$}. \mbox{$\vartheta\in[0,\pi]$} is the polar (colatitudinal) coordinate, and
$\varphi\in[0,2\pi)$ -- the azimuthal (longitudinal) coordinate.

The scalar product of spherical functions~$\Phi$, $\Psi$ is defined by
$$
\left<\Phi,\Psi\right>=\int_{\mathcal S}\overline{\Phi(x)}\,\Psi(x)\,d\sigma(x),
$$
and $\mathcal L^2(\mathcal S)$ is the set of functions~$\Psi$ such that
$$
\|\Psi\|_{\mathcal L^2(\mathcal S)}:=\left<\Psi,\Psi\right>^{\frac{1}{2}}<\infty.
$$
Note the difference to the notation usually used in research concerning $n$-dimensional spheres, where the integral in the definition of the scalar product is divided by~$\Sigma_n$, the Lebesgue measure of~$\mathcal S^n$, see e.g. \cite{IIN15CWT}.
The spherical harmonics $Y_l^k$, $l\in\mathbb N_0$, $k=-l,\mbox{$-l+1$},\dots,l$, are given by
   \begin{align*}
   Y_l^k(\vartheta,\varphi)&=(-1)^k\sqrt{\frac{2l+1}{4\pi}\frac{(l-k)!}{(l+k)!}}\,P_l^k(\cos\vartheta)\,e^{ik\varphi},\\
   Y_l^{-k}&=\overline{Y_l^k}
   \end{align*}
for $0\leq k\leq  l$, where $P_l^k$ are the (associated) Legendre polynomials (functions),
   \begin{align*}
   P_l^{\,0}(t)&=P_l(t)=\frac{1}{2^l l!}\frac{d^l}{dt^l}(t^2-1)^l&&\text{(Legendre polynomials)}\\
   P_l^k(t)&=(-1)^k(1-t^2)^{k/2}\frac{d^k}{dt^k}P_l(t)\quad\text{for }0<k\leq l&&\text{(associated Legendre functions).}
   \end{align*}
The spherical harmonics build an orthonormal basis for $\mathcal L^2(\mathcal S)$. For an $\mathcal L^2(\mathcal S)$--function~$\Psi$, the series
   \[
   \sum_{l=-\infty}^\infty\sum_{k=-l}^{l}\left<Y_l^k,\Psi\right>Y_l^k
   \]
is the Fourier series of~$\Psi$ in terms of the spherical harmonics. The constants $\widehat\Psi_l^k=\left<Y_l^k,\Psi\right>$ are called the Fourier coefficients of~$\Psi$.

$\mathbb T$ denotes the one-dimensional torus $\mathbb T=[0,2\pi)$. Functions on~$\mathbb T$ will be identified with their $2\pi$-periodic extension on~$\mathbb R$, such that $\mathcal C(\mathbb T)$ denotes the class of continuous functions on~$\mathbb T$ that can be continuously extended. An $\mathcal L^1(\mathbb T)$-function can be represented as the Fourier series
$$
f(t)\sim\frac{1}{2\pi}\sum_{k=-\infty}^\infty\widehat f_k\,e^{ikt},
$$
where the Fourier coefficients~$\widehat f_k$, $k\in\mathbb Z$, are given by
$$
\widehat f_k=\int_0^{2\pi}f(t)\,e^{-ikt}\,dt.
$$
Similarly, the Fourier transform of an $\mathcal L^1(\mathbb R)$ function~$F$ is defined by
$$
\widehat F(\omega)=\int_{-\infty}^\infty F(t)\,e^{-i\omega t}\,dt,\qquad\sigma\in\mathbb R,
$$
and it can be inverted by
$$
F(t)=\frac{1}{2\pi}\int_{-\infty}^\infty\widehat F(\omega)\,e^{i\omega t}\,d\sigma.
$$
The functions we consider in the present paper are both integrable and continuous such that integrals and series in the Fourier and the inverse Fourier transform are (pointwise) convergent. An $\mathcal L^1(\mathbb R)$-function~$F$ can be periodized to an $\mathcal L^1(\mathbb T)$-function~$f$ by
$$
f(t)=\sum_{j=-\infty}^\infty f(t+2j\pi),\qquad t\in\mathbb T,
$$
and their Fourier transforms are linked by the following Poisson summation formula.

\begin{thm}\label{thm:PSF}
Let $F\in\mathcal L^1(\mathbb R)$. Then its periodization~$f$ exists for almost every $t\in\mathbb T$, it is integrable and $\|f\|_{\mathcal L^1(\mathbb T)}\leq\|F\|_{\mathcal L^1(\mathbb R)}$. Further,
$$
\widehat f_k=\widehat F(k).
$$
\end{thm}

A family of vectors $\{g_\iota,\,\iota\in I\}\subset\mathcal{H}$ in a Hilbert space $\mathcal{H}$ indexed by a measure space $I$ with a positive measure $\mu$ is called a~\emph{frame with weight}~$\mu$ if the mapping $\iota\mapsto g_\iota$ is weakly measurable, i.e. $\iota\mapsto\left<g_\iota,u\right>$ is measurable and
$$
A\,\|u\|^2\leq \int_{I}|\left<g_\iota,u\right>|^2\,d\mu(x) \leq B\,\|u\|^2
$$
for all $u\in\mathcal{H}$ and $0<A\leq B$. If $\{g_\iota,\,\iota\in I\}$ is a frame, then the mapping
$$
u\mapsto\{\left<g_\iota,u\right>,\,\iota\in I\}
$$
is invertible.

The wavelet definition we use originates from~\cite{IIN16FDW}. In the case of the two-dimensional sphere and the weight function $\rho\mapsto\frac{1}{\rho}$, it is the following one.

\begin{df}\label{def:wavelet} The family $\{\Psi_\rho\}_{\rho\in\mathbb R_+}\subseteq\mathcal L^2(\mathbb S^2)$ is called a wavelet (family) of order~$m$ if it satisfies
\begin{equation}\label{eq:wavelet_condition}
A(2l+1)\leq\sum_{k=-l}^l\int_0^\infty\left|(\widehat\Psi_\rho)_l^k\right|^2\,\frac{d\rho}{\rho}\leq B(2l+1)
\end{equation}
for some positive constants~$A$ and~$B$ independent of $l\in\mathbb{N}_0$, $l>m$, and
\begin{equation}\label{eq:wavelet_0_condition}
\sum_{k=-l}^l\int_0^\infty\left|(\widehat\Psi_\rho)_l^k\right|^2\,\frac{d\rho}{\rho}=0
\end{equation}
for $l\in\mathbb N_0$, $l\leq m$.
\end{df}

The wavelet transform of an $\mathcal L^2(\mathcal S^2)$-function is given by
$$
\mathcal W_\Psi f(\rho,g)=\frac{1}{4\pi}\int_{\mathcal S^n}\overline{\Psi_\rho(g^{-1}x)}\,f(x)\,d\sigma_n(x),\qquad g\in SO(3),
$$
and it is invertible by methods known from the frame theory. Moreover, only a discrete set of the wavelet coefficients~$\mathcal W_\Psi f(\rho,g)$ is necessary to reconstruct the analyzed signal. The first step is a scale discretization. According to \cite[Theorem~3.1]{IIN16FDW} there exist constants~$\mathfrak a_0$ and~$X$ such that for any sequence $\mathcal R=(\rho_j)_{j\in\mathbb N_0}$ with $\rho_0\geq\mathfrak a_0$ and $1<\rho_j/\rho_{j+1}<X$ the family $\{\Psi_{\rho_j}(g^{-1}\circ),\,\rho_j\in\mathcal R,\,g\in SO(3)\}$ is a frame for~$\mathcal L^2(\mathcal S^2)$. In order to characterize discrete sets of rotations, we need the following definition (compare \cite[Definition~4.1]{IIN16FDW}).

\begin{df}\label{df:grid}$\Lambda$ is a grid of type $(\delta_2,\delta_1)$ if it is a discrete measurable set of rotations in~$SO(3)=\mathcal S^2\times\mathcal S^1$, constructed in the following way. Let $\mathcal P_2=\{\mathcal O_{\alpha_2}^2:\,\alpha_2=1,\dots,K_2\}$ be a measurable partition of~$\mathcal S^2$ into simply connected sets such that the diameter of each set (measured in the geodesic distance) is not larger than~$\delta_2$. Choose from each of the sets~$\mathcal O_{\alpha_2}^2$ an arbitrary point~$x_{\alpha_2}^2=(\vartheta_{\alpha_2}^2,\varphi_{\alpha_2}^2)$. Now, for a fxed~$\alpha_2$, let $\mathcal P_1(\alpha_2)=\{\mathcal O_{(\alpha_2,\alpha_1)}^1:\,\alpha_1=1,\dots,K_1(\alpha_2)\}$ be a  measurable partition of~$\mathcal S^1$ into $K_1(\alpha_2)$ simply connected sets of a diameter not larger than~$\delta_1$. Choose from each of the sets~$\mathcal O_{(\alpha_2,\alpha_1)}$, $\alpha_1=1,\dots,K_1(\alpha_2)$, an arbitrary point~$\varphi_{(\alpha_2,\alpha_1)}^1$. Then, $\Lambda$ is the set of rotations given by
$$
g_{(\alpha_2,\alpha_1)}=g_1(\varphi_{(\alpha_2,\alpha_1)}^1)\,g_1(\vartheta_{\alpha_2}^2)\,g_2(\varphi_{\alpha_2}^2),
$$
where~$g_j(\beta)$ is the rotation in the plane $(\xi_j,\xi_{j+1})$ with the rotation angle~$\beta$.
\end{df}

\begin{figure}\label{fig:grid_Lambda}
\vspace{-8em}
\includegraphics[angle=270,scale=0.6]{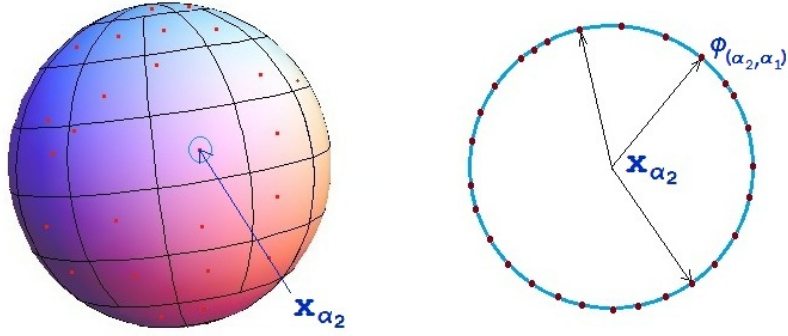}\vspace{-10em}
\caption{Grid construction}
\end{figure}

Loosely speaking, one has a tight enough discrete set of points $\{x_{\alpha_2}\}\in\mathcal S^2$, and for each point a tight enough discrete set of rotations around it. $\Lambda$ is a set of $SO(3)$-rotations being a composition of a tranlastion to a point~$x_{\alpha_2}$ and a rotation around this point, as illustrated in  Figure~\ref{fig:grid_Lambda}.

The next statement reflects the content of \cite[Theorem~4.2]{IIN16FDW} for the case of the two-dimensional sphere.

\begin{thm} Let $\Psi_\rho$ be a $\mathcal C^1$--wavelet family with the property that $\{\Psi_{\rho_j,x},\,j\in\mathbb N_0,\,x\in\mathcal S^2\}$ is a semi--continuous frame. Then, for each $j\in\mathbb N_0$ there exist numbers $\delta_2^j$, $\delta_1^j$ such that
$$
\{\Psi_{\rho_j}(g_{(\alpha_2^j,\alpha_1^j)}^{-1}\circ),\,j\in\mathbb N_0,\,g_{(\alpha_2^j,\alpha_1^j)}\in\Lambda^j\}
$$
is a frame for~$\mathcal L^2(\mathcal S^n)$, provided that $\Lambda^j$ is a grid of type $(\delta_2^j,\delta_1^j)$.
\end{thm}

\section{Construction of the wavelet}\label{sec:construction}

In this section, we shall present a wavelet family with separated variables and a steerable angular selectivity. The construction is based on the Poisson kernel,
\begin{equation}\label{eq:Poisson_kernel}
p_\rho(x)=\frac{1}{4\pi}\frac{1-r^2}{(1-2r\cos\vartheta+r^2)^{3/2}}=\frac{1}{4\pi}\sum_{l=0}^\infty(2l+1)r^l P_l(\cos\vartheta),
\end{equation}
where $r=e^{-\rho}$, since our aim is to obtain a wavelet with a simple representation as a function of the spherical variables. On the other hand, in order to prove the wavelet property, some knowledge about the Fourier coefficients is necessary.

\begin{figure}
\centering
\centering\vspace{-8em}
\includegraphics[angle=270,scale=0.4]{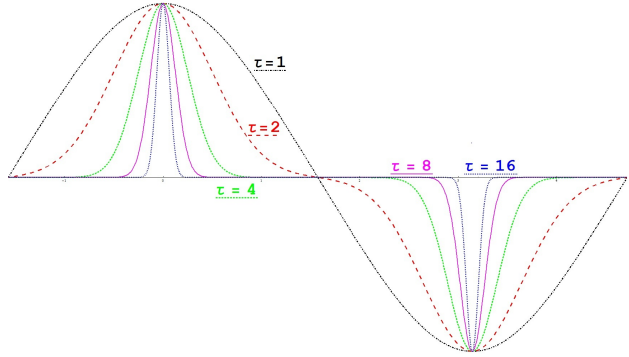}\vspace{-6em}
\caption{$f_\tau$ for $\tau=1,\,2,\,4,\,8,\,16$, $\varphi\in[-\pi/2,3\pi/2]$}\label{fig:f_tau}
\end{figure}

The idea is to choose a periodized difference of Gaussians dependent on a parameter $\tau\in[1,\infty)$ as the function of the longitudinal variable~$\varphi$,
\begin{align*}
F_\tau(\phi)&=e^{-\frac{\tau^2\phi^2}{2}}-e^{-\frac{\tau^2(\phi-\pi)^2}{2}},\qquad\phi\in\mathbb R,\\
f_\tau(\varphi)&=\sum_{j\in\mathbb Z}F_\tau(\varphi+2j\pi),\qquad\varphi\in[0,2\pi).
\end{align*}
Figure~\ref{fig:f_tau} presents the function~$f_\tau$ for several values of~$\tau$.

The challenge is to find an appropriate function of the colatitudinal variable~$\vartheta$ such that the inequalities~\eqref{eq:wavelet_condition} are satisfied. The result of our investigation is that the following families:
\begin{align*}
\Omega(x)&=\Omega_\rho^\tau(\vartheta,\varphi)=\omega_\rho(\vartheta)\cdot f_\tau(\varphi),\qquad
   \omega_\rho(\vartheta)=\rho\cdot\sin^5\vartheta\cdot r\cdot\frac{\partial}{\partial r}\left[r\cdot\frac{\partial}{\partial r}\left[p_\rho(x)\right]\right],\\
\Upsilon(x)&=\Upsilon_\rho^\tau(\vartheta,\varphi)=\upsilon_\rho(\vartheta)\cdot f_\tau(\varphi),\qquad
   \upsilon_\rho(\vartheta)=\rho\cdot\sin^5\vartheta\cdot r^2\cdot\frac{\partial^2}{\partial r^2}\left[p_\rho(x)\right],
\end{align*}
$r=e^{-\rho}$, are wavelets.

In order to prove the wavelet condition~\eqref{eq:wavelet_condition}, we need an estimation of the Fourier coefficients (with respect to the spherical harmonics) of~$\Omega$ and $\Upsilon$, and since they are functions with separated variables, we first consider the Fourier series expansion (in exponential functions) of~$f_\tau$. Since
$$
\widehat F_\tau(\omega)=\frac{\sqrt{2\pi}}{\tau}\left(e^{-\frac{\omega^2}{2\tau^2}}-e^{-\frac{\omega^2}{2\tau^2}-i\omega\pi}\right),
$$
by the Poisson summation formula the Fourier coefficients of~$f_\tau$ are given by
$$
(\widehat f_\tau)_k=\begin{cases}0,&\text{for even }k,\\\frac{2\sqrt{2\pi}}{\tau}\,e^{-\frac{k^2}{2\tau^2}}&\text{for odd }k.\end{cases}
$$
Thus,
\begin{equation}\label{eq:f_alpha_series}
f_\tau(\varphi)=\frac{1}{\tau}\sqrt\frac{2}{\pi}\sum_{k\in2\mathbb Z+1}e^{-\frac{k^2}{2\tau^2}}\,e^{ik\varphi}.
\end{equation}
In order to estimate the Fourier coefficients of~$\Omega$ and~$\Upsilon$, we estimate the magnitude of the coefficients of
\begin{align}
\omega_\rho(\vartheta)&=\frac{\rho\sin^5\vartheta}{4\pi}\sum_{l=1}^\infty(2l+1)l^2r^lP_l(\cos\vartheta)\label{eq:omega_series}\\
&=-\frac{\rho r\left[r(10-19r^2+r^4)-(3-14r^2+5r^4)\cos\vartheta-r(9-r^2)\cos^2\vartheta\right]\sin^5\vartheta}{4\pi(1-2r\cos\vartheta+r^2)^{7/2}},\label{eq:omega_rational}
\end{align}
and
\begin{align}
\upsilon_\rho(\vartheta)&=\frac{\rho\sin^5\vartheta}{4\pi}\sum_{l=2}^\infty(2l+1)l(l-1)r^lP_l(\cos\vartheta)\label{eq:upsilon_series}\\
&=-\frac{\rho r^2\left[5-23r^2+2r^4+4r(7+r^2)\cos\vartheta-(15+r^2)\cos^2\vartheta\right]\sin^5\vartheta}{4\pi(1-2r\cos\vartheta+r^2)^{7/2}},\label{eq:upsilon_rational}
\end{align}
$r=e^{-\rho}$, developed into series with respect to the associated Legendre polynomials of an odd order~$k$.

\begin{lem}\label{lem:omega_series_P1} For $r\in(0,1)$,
\begin{align*}
\omega_\rho(\vartheta)&=\frac{\rho}{4\pi}\sum_{l=1}^\infty\beta_l P_l^1(\cos\vartheta)\quad\text{and}\\
\upsilon_\rho(\vartheta)&=\frac{\rho}{4\pi}\sum_{l=2}^\infty\gamma_l P_l^1(\cos\vartheta)
\end{align*}
with
\begin{equation}\begin{split}\label{eq:omega_seriesP1}
\beta_l&=\frac{(l+2)(l+3)(l+4)(l+5)^3\,r^{l+5}}{(2l+3)(2l+5)(2l+7)(2l+9)}-\frac{(l+2)(l+3)^3(5l^2+21l-8)\,r^{l+3}}{(2l-1)(2l+3)(2l+5)(2l+9)}\\
&+\frac{2(l+1)^2(10l^5+65l^4+83l^3-134l^2-186l+36)\,r^{l+1}}{(2l-3)(2l-1)(2l+3)(2l+5)(2l+7)}\\
&-\frac{2(l-1)^2(5l^4-2l^3-41l^2+14l+60)\,r^{l-1}}{(2l-5)(2l-3)(2l+3)(2l+5)}\\
&+\frac{(l-3)^2(l-2)(l-1)(20l^4-100l^3+39l^2+239l-48)\,r^{l-3}}{(2l-7)(2l-5)(2l-3)(2l-1)(1+2l)(2l+3)}\\
&-\frac{(l-5)^2(l-4)(l-3)(l-2)(l-1)\,r^{l-5}}{(2l-7)(2l-5)(2l-3)(2l-1)}
\end{split}\end{equation}
and
\begin{equation}\begin{split}\label{eq:upsilon_seriesP1}
\gamma_l&=\frac{(l+2)(l+3)(l+4)^2(l+5)^2\,r^{l+5}}{(2l+3)(2l+5)(2l+7)(2l+9)}-\frac{(l+2)^2(l+3)^2(5l^2+21l-8)\,r^{l+3}}{(2l-1)(2l+3)(2l+5)(2l+9)}\\
&+\frac{2l(l+1)(10l^5+65l^4+83l^3-134l^2-186l+36)\,r^{l+1}}{(2l-3)(2l-1)(2l+3)(2l+5)(2l+7)}\\
&-\frac{2(l-2)(l-1)(5l^4-2l^3-41l^2+14l+60)\,r^{l-1}}{(2l-5)(2l-3)(2l+3)(2l+5)}\\
&+\frac{(l-4)(l-3)(l-2)(l-1)(20l^4-100l^3+39l^2+239l-48)\,r^{l-3}}{(2l-7)(2l-5)(2l-3)(2l-1)(1+2l)(2l+3)}\\
&-\frac{(l-6)(l-5)(l-4)(l-3)(l-2)(l-1)\,r^{l-5}}{(2l-7)(2l-5)(2l-3)(2l-1)}
\end{split}\end{equation}
for $l\geq9$. The coefficients $\beta_l$, $l=1,\,2,\dots,\,8$, and $\gamma_l$, $l=2,\,3,\dots,\,8$, do not vanish.
\end{lem}

\begin{bfseries}Proof. \end{bfseries} By the repeated application of \cite[formulae~(8.733.4) and~(8.735.5)]{GR}:
\begin{align*}
\sin\vartheta P_l^{k-1}(\cos\vartheta)&=\frac{P_{l-1}^k(\cos\vartheta)-P_{l+1}^k(\cos\vartheta)}{2l+1},\\
\sin\vartheta P_l^{k+1}(\cos\vartheta)&=\frac{(l-k)(l-k+1)}{2l+1}P_{l+1}^k(\cos\vartheta)-\frac{(l+k)(l+k+1)}{2l+1}P_{l-1}^k(\cos\vartheta)
\end{align*}
we obtain
\begin{align*}
(2l&+1)\sin^5\vartheta P_l^0(\cos\vartheta)
   =\frac{(l-3)(l-2)(l-1)l}{(2l-7)(2l-5)(2l-3)(2l-1)}P_{l-5}^1(\cos\vartheta)\\
&-\frac{(l-1)l(5l^2-9l-26)}{(2l-7)(2l-3)(2l-1)(2l+3)}P_{l-3}^1(\cos\vartheta)\\
&+\frac{2(10l^5+15l^4-77l^3-93l^2+121l+60)}{(2l-5)(2l-3)(2l+1)(2l+3)(2l+5)}P_{l-1}^1(\cos\vartheta)\\
&-\frac{2(5l^4+18l^3-17l^2-54l+36)}{(2l-3)(2l-1)(2l+5)(2l+7)}P_{l+1}^1(\cos\vartheta)\\
&+\frac{(l+1)(l+2)(20l^4+140l^3+219l^2-67l-60)}{(2l-1)(2l+1)(2l+3)(2l+5)(2l+7)(2l+9)}P_{l+3}^1(\cos\vartheta)\\
&-\frac{(l+1)(l+2)(l+3)(l+4)}{(2l+3)(2l+5)(2l+7)(2l+9)}P_{l+5}^1(\cos\vartheta)
\end{align*}
for $l\geq4$. A substitution of this formula to~\eqref{eq:omega_series} resp.~\eqref{eq:upsilon_series} and index shifts yield coefficients of~$\omega_\rho$ resp.~$\upsilon_\rho$ as in~\eqref{eq:omega_seriesP1} resp.~\eqref{eq:upsilon_seriesP1} for $l\geq9$. The coefficients~$\beta_l$, $\gamma_l$, $l=1,\,2,\dots,\,8$, are computed in the same manner. The difference with respect to the case $l\geq9$ is caused by vanishing polynomials $P_0^1$.\hfill$\Box$\\

\begin{lem}\label{lem:lower_bound}The Fourier coefficients~$\widehat\Omega_l^1$, $l=1,2,\dots$, and $\widehat\Upsilon_l^1$, $l=2,3,\dots$, satisfy
\begin{equation}\label{eq:lower_bound}
(2l+1)A_\Omega\leq\int_0^\infty\left|\widehat\Omega_l^1\right|^2\frac{d\rho}{\rho},\qquad(2l+1)A_\Upsilon\leq\int_0^\infty\left|\widehat\Upsilon_l^1\right|^2\frac{d\rho}{\rho}
\end{equation}
for some constants~$A_\Omega$, $A_\Upsilon$ independent of~$l$.
\end{lem}

\begin{bfseries}Proof.\end{bfseries} It follows from Lemma~\ref{lem:omega_series_P1} and formula~\eqref{eq:f_alpha_series} that
\begin{equation*}\begin{split}
\Omega(x)&=\omega_\rho(\vartheta)\cdot f_\tau(\varphi)=\frac{\rho}{4\pi}\sum_{l=1}^\infty\beta_lP_l^1(\cos\vartheta)
   \cdot\frac{1}{\tau}\sqrt\frac{2}{\pi}\sum_{k\in2\mathbb Z+1}e^{-\frac{k^2}{2\tau^2}}\,e^{ik\varphi}\\
&=\frac{\rho}{4\tau\pi}\sqrt\frac{2}{\pi}\sum_{l,k}\beta_l\,e^{-\frac{k^2}{2\tau^2}}\,P_l^1(\cos\vartheta)\,e^{ik\varphi}.
\end{split}\end{equation*}
Since
$$
P_l^1(\cos\vartheta)\,e^{i\varphi}=-\sqrt\frac{4\pi\cdot l(l+1)}{2l+1}\,Y_l^1(\vartheta,\varphi)
$$
and the spherical harmonics are orthogonal to one another,
$$
\widehat\Omega_l^1=-\frac{\rho}{4\tau\pi}\sqrt\frac{2}{\pi}\cdot\beta_l\,e^{-\frac{1}{2\tau^2}}\cdot\sqrt\frac{4\pi\cdot l(l+1)}{2l+1}
   =-\frac{\rho}{\tau\pi}\sqrt\frac{l(l+1)}{2(2l+1)}\cdot\beta_l\,e^{-\frac{1}{2\tau^2}}.
$$
Similarly,
$$
\widehat\Upsilon_l^1=-\frac{\rho}{\tau\pi}\sqrt\frac{l(l+1)}{2(2l+1)}\cdot\gamma_l\,e^{-\frac{1}{2\tau^2}}.
$$
For $l\geq9$, the integrals
$$
\frac{1}{2l+1}\int_0^\infty\left|\widehat\Omega_l^1\right|^2\frac{d\rho}{\rho}\qquad\text{and}\qquad\frac{1}{2l+1}\int_0^\infty\left|\widehat\Upsilon_l^1\right|^2\frac{d\rho}{\rho}
$$
are positive rational functions of~$l$. It can be verified that
$$
\int_0^\infty\rho\,\beta_l^2\,d\rho\qquad\text{and}\qquad\int_0^\infty\rho\,\gamma_l^2\,d\rho
$$
are continuous for $l\geq9$ and that their limits for $l\to\infty$ both equal~$\frac{1}{256}$. Consequently,
$$
A_\Omega^9:=\min_{l\in[9,\infty)}\frac{1}{2l+1}\int_0^\infty\left|\widehat\Omega_l^1\right|^2\frac{d\rho}{\rho}\qquad\text{and}\qquad
A_\Upsilon^9:=\min_{l\in[9,\infty)}\frac{1}{2l+1}\int_0^\infty\left|\widehat\Upsilon_l^1\right|^2\frac{d\rho}{\rho}
$$
both exist and are positive. Similarly, for $l=1,\,2,\dots,\,8$ the numbers
$$
A_\Omega^l:=\frac{1}{2l+1}\int_0^\infty\left|\widehat\Omega_l^1\right|^2\frac{d\rho}{\rho}
$$
and for $l=2,\,3,\dots,\,8$ the numbers
$$
A_\Upsilon^l:=\frac{1}{2l+1}\int_0^\infty\left|\widehat\Upsilon_l^1\right|^2\frac{d\rho}{\rho}
$$
are positive. Inequalities~\eqref{eq:lower_bound} are thus satisfied with $A_\Omega=\min\{A_\Omega^1,\,A_\Omega^2,\dots,A_\Omega^8,\,A_\Omega^9\}$ and $A_\Upsilon=\min\{A_\Upsilon^2,\,A_\Upsilon^3,\dots,\,A_\Upsilon^8,\,A_\Upsilon^9\}$.\hfill$\Box$\\

\begin{lem}\label{lem:psi_Plk_integral}Let $k\in\mathbb N$. Then
\begin{align}
&\left|\int_0^\pi\omega_\rho(\vartheta)P_l^k(\cos\vartheta)\sin\vartheta\,d\vartheta\right|\leq\frac{3\rho r}{2\sqrt{1-r^2}}\sqrt\frac{(l+k)!}{k(l-k)!}
   \qquad\text{and}\label{eq:upper_bound_omega_Plk}\\
&\left|\int_0^\pi\upsilon_\rho(\vartheta)P_l^k(\cos\vartheta)\sin\vartheta\,d\vartheta\right|\leq\frac{3\rho r^2}{\sqrt{1-r^2}}\sqrt\frac{(l+k)!}{k(l-k)!}
   \label{eq:upper_bound_upsilon_Plk}
\end{align}
\end{lem}

\begin{bfseries}Proof. \end{bfseries}Denote the left-hand-sides of~\eqref{eq:upper_bound_omega_Plk} and~\eqref{eq:upper_bound_upsilon_Plk} by~$\widetilde\omega_l^k$ and~$\widetilde\upsilon_l^k$, respectively.
They are bounded from above by
\begin{align*}
&\widetilde\omega_l^k\leq\frac{\rho r}{4\pi}\int_0^\pi\frac{\sin^6\vartheta\,|p_\omega(r,\vartheta)|\left|P_l^k(\cos\vartheta)\right|}{(1-2r\cos\vartheta+r^2)^{7/2}}\,d\vartheta
   \qquad\text{and}\\
&\widetilde\upsilon_l^k
   \leq\frac{\rho r^2}{4\pi}\int_0^\pi\frac{\sin^6\vartheta\,|p_\upsilon(r,\vartheta)|\left|P_l^k(\cos\vartheta)\right|}{(1-2r\cos\vartheta+r^2)^{7/2}}\,d\vartheta,
\end{align*}
where
\begin{align*}
&p_\omega(r,\vartheta)=r(10-19r^2+r^4)-(3-14r^2+5r^4)\cos\vartheta-r(9-r^2)\cos^2\vartheta\qquad\text{and}\\
&p_\upsilon(r,\vartheta)=5-23r^2+2r^4+4r(7+r^2)\cos\vartheta-(15+r^2)\cos^2\vartheta,
\end{align*}
see \eqref{eq:omega_rational} and~\eqref{eq:upsilon_rational}. The H\"older inequality together with
$$
\int_0^\pi\frac{[P_l^k(\cos\vartheta)]^2}{\sin\vartheta}\,d\vartheta=\frac{(l+k)!}{k(l-k)!}
$$
(see~\cite[formula~(7.122.1)]{GR}, note that $2$ in the denominator disappears because of the extended integration interval) yield
\begin{align*}
\widetilde\omega_l^k\leq&\frac{\rho r}{4\pi}\sup_{\genfrac{}{}{0pt}{}{r\in(0,1)}{\vartheta\in[0,\pi]}}\left(\sqrt{\sin\vartheta}\,|p_\omega(r,\vartheta)|\right)\\
&\cdot\left[\int_0^\pi\frac{\sin^{12}\vartheta}{(1-2r\cos\vartheta+r^2)^7}\,d\vartheta\right]^{1/2}
   \cdot\left[\int_0^\pi\frac{\left[P_l^k(\cos\vartheta)\right]^2}{\sin\vartheta}\,d\vartheta\right]^{1/2}\\
\leq&\frac{\rho r}{4\pi}\cdot21\cdot\left[\frac{231\pi}{1024(1-r^2)}\cdot\frac{(l+k)!}{k(l-k)!}\right]^{1/2}
   \leq\frac{3\rho r}{2\sqrt{(1-r^2)}}\sqrt\frac{(l+k)!}{k(l-k)!}.
\end{align*}
Estimation~\eqref{eq:upper_bound_upsilon_Plk} is obtained in a similar way. It is to be respected that the function $(r,\vartheta)\mapsto\sqrt{\sin\vartheta}\,|p_\upsilon(r,\vartheta)|$ is bounded from above by~$41$.\hfill$\Box$\\

\begin{thm}The function family~$\{\Omega_\rho^\tau\}$ is a wavelet of order~$0$, and the function family~$\{\Upsilon_\rho^\tau\}$ is a wavelet of order~$1$ according to Definition~\ref{def:wavelet}.
\end{thm}

\begin{bfseries}Proof. \end{bfseries}The lower bound inequality in~\eqref{eq:wavelet_condition} is satisfied by Lemma~\ref{lem:lower_bound} for \mbox{$l\geq1$} in the case of the wavelet~$\{\Omega_\rho^\tau\}$, respectively for $l\geq2$ in the case of the wavelet~$\{\Upsilon_\rho^\tau\}$. In order to obtain an estimation from above, note that for odd $k$'s
\begin{align*}
\left|\Omega_l^k\right|&=\sqrt{\frac{2l+1}{4\pi}\frac{(l-|k|)!}{(l+|k|)!}}
   \cdot\int_0^\pi\omega_\rho(\vartheta)\,P_l^{|k|}(\cos\vartheta)\sin\vartheta\,d\vartheta\cdot\int_0^{2\pi}f_\tau(\varphi)\,e^{-ik\varphi}\,d\varphi\\
&\leq\sqrt{\frac{2l+1}{4\pi}\frac{(l-|k|)!}{(l+|k|)!}}\cdot\frac{3\rho r}{2\sqrt{1-r^2}}\sqrt\frac{(l+|k|)!}{|k|(l-|k|)!}\cdot\frac{2\sqrt{2\pi}}{\tau}\,e^{-\frac{k^2}{2\tau^2}}\\
&=\frac{3\rho r}{\tau}\sqrt\frac{2l+1}{2|k|(1-r^2)}\,e^{-\frac{k^2}{2\tau^2}}
\end{align*}
and, analogously,
$$
\left|\Upsilon_l^k\right|\leq\frac{6\rho r^2}{\tau}\sqrt\frac{2l+1}{2|k|(1-r^2)}\,e^{-\frac{k^2}{2\tau^2}}.
$$
Thus, since
$$
\int_0^\infty\frac{\rho\,e^{-2\rho}}{1-e^{-2\rho}}\,d\rho=\frac{\pi^2}{24}\qquad\text{and}\qquad
    \int_0^\infty\frac{\rho\,e^{-4\rho}}{1-e^{-2\rho}}\,d\rho=\frac{\pi^2-6}{24},
$$
we obtain for each~$l$
$$
\sum_{k=-l}^l\left|\Omega_l^k\right|^2\frac{d\rho}{\rho}\leq\frac{9\pi^2(2l+1)}{2\cdot24\,\tau^2}\sum_{k\in2\mathbb Z+1}\frac{e^{-\frac{k^2}{\tau^2}}}{|k|}
   \leq\frac{2(2l+1)}{\tau^2}\sum_{k\in2\mathbb Z+1}\frac{e^{-\frac{k^2}{\tau^2}}}{|k|}
$$
and
$$
\sum_{k=-l}^l\left|\Upsilon_l^k\right|^2\frac{d\rho}{\rho}\leq\frac{36(\pi^2-6)(2l+1)}{2\cdot24\,\tau^2}\sum_{k\in2\mathbb Z+1}\frac{e^{-\frac{k^2}{\tau^2}}}{|k|}
   \leq\frac{3(2l+1)}{\tau^2}\sum_{k\in2\mathbb Z+1}\frac{e^{-\frac{k^2}{\tau^2}}}{|k|}.
$$
Further, since
$$
t\mapsto\frac{e^{-\frac{t^2}{\tau^2}}}{t}
$$
is a decreasing function for $t\in(0,\infty)$, the estimation
$$
\sum_{k\in2\mathbb N_0+1}\frac{e^{-\frac{k^2}{\tau^2}}}{|k|}\leq\frac{1}{2}\int_1^\infty\frac{e^{-\frac{t^2}{\tau^2}}}{t}\,dt
$$
is valid, and since
$$
\int_1^\infty\frac{e^{-\frac{t^2}{\tau^2}}}{t}\,dt=\int_\frac{1}{\tau}^\infty\frac{e^{-u^2}}{u}\,du\leq\int_\frac{1}{\tau}^1\frac{du}{u}+\int_0^\infty e^{-u^2}du
   =\ln\tau+\frac{\sqrt\pi}{2},
$$
the upper bound inequality in~\eqref{eq:wavelet_condition} holds with
$$
B_\Omega=\frac{2}{\tau^2}\left(\ln\tau+\frac{\sqrt\pi}{2}\right)\qquad\text{and}\qquad B_\Upsilon=\frac{3}{\tau^2}\left(\ln\tau+\frac{\sqrt\pi}{2}\right)
$$
for the families $\{\Omega_\rho^\tau\}$ and $\{\Upsilon_\rho^\tau\}$, respectively. Condition~\eqref{eq:wavelet_0_condition} follows from the representations~\eqref{eq:omega_series} and~\eqref{eq:upsilon_series}. \hfill$\Box$

\section{Multiselectivity scheme}\label{sec:multiselectivity}

Our goal is to define an angular multiselectivity analysis similar to that presented in~\cite{AJ03book,AJ03,JA07}. The main difference between our work and that of Antoine and Jacques is that in the $\mathbb R^2$-case, an exact reconstruction is possible with the use of biorthogonal wavelets and scaling functions, whereas the spherical wavelet transform is invertible by frame methods.

A multiselectivity analysis is defined a priori for the discrete case. The continuous wavelet transform is performed with respect to a fixed wavelet, and it is certainly possible to choose a function with a proper resolution power, but the choice must be based on some global estimations. In the case of a wavelet transform for a discrete set of scales and positions, the frame property remains unaltered if a wavelet is chosen for each of the parameters separately. The angular selectivity can be adapted to the features detected in each point of the analyzed signal by varying the parameter~$\tau$. However, the discretization density of the sphere~$\mathcal S^2$ must be chosen in such a way that the frame property is preserved in the case that a wavelet with the highest angular resolution is used.

More exactly, according to the proof of~\cite[Theorem~4.2]{IIN16FDW}, the error~$\mathcal E^2$ made by the discretization of the set of wavelet coefficients is bounded by
$$
|\mathcal E^2|\leq4\pi E_1+E_2,
$$
with
$$
E_J=\mathfrak c\cdot\|\Psi_{\rho_j}\|_\infty\cdot\|\nabla^\ast\Psi_{\rho_j}\|_\infty\cdot\delta_J^j\cdot\|f\|^2,\qquad J=1,2,
$$
where~$\Psi_{\rho_j}$ is the wavelet on the scale~$\rho_j$, and $\delta_J^j$ is an upper bound for the diameter of the discretization sets on the scale~$\rho_j$, and $f$ denotes the analyzed signal. Since in general $|\nabla^\ast\Psi_{\rho_j}|$ increases with increasing~$\tau$, and $|\mathcal E^2|\leq2^{-j-1}\delta\|f\|_2^2$ is required, an upper bound~$\widetilde\tau$ for~$\tau$ must be set, and $\delta_2^j$ computed with respect to~$\Psi_\rho^{\widetilde\tau}$.

Suppose the discretization of scales and $\mathcal S^2$-positions is performed. For each scale $\rho_j$ and each position~$x_{\alpha_2}^2$, the 'best frame' for~$f$, i.e. the best angular selectivity level is to be chosen (see~\cite[Section~5]{AJ03}). The best match between the image and the wavelet can be quantified by
$$
\tau(\rho_j,x_{\alpha_2}^2)=\underset{\tau\in[1,\widetilde\tau]}{\text{arg max}}\,\max_{\genfrac{}{}{0pt}{}{g\in SO(3)}{g^{-1}\widehat e=x_{\alpha_2}^2}}
\frac{\left|\left<\Psi_{\rho_j}^\tau(g^{-1}\circ),f\right>\right|}{\|\Psi_{\rho_j}\|^2},
$$
where~$\widehat e$ denotes the north pole of the sphere. Alternatively, one can a priori choose a discrete (and bounded from above) set~$T$ of possible~$\tau$'s, and replace
$$
\underset{\tau\in[1,\widetilde\tau]}{\text{arg max}}
$$
in the above formula by
$$
\underset{\tau\in T}{\text{arg max}}.
$$
Note that according to the proof of \cite[Theorem~4.2]{IIN16FDW}, the discretization error~$E_1$ is an upper bound for a weighted sum of the discretization errors~$e_{\alpha_1}^1(\alpha_2)$ made in each grid point $\alpha_2\in\mathcal S^2$ that can be estimated by
$$
e_{\alpha_1}^1(\alpha_2)\leq\mathfrak c\cdot\left\|\sup|\Psi_{\rho_j}|\right\|_2\cdot\left\|\sup|\nabla^\ast\Psi_{\rho_j}|\right\|_2\cdot\|f\|_2^2\cdot\delta_1^j,
$$
(absolute values are missing in~\cite[formula~(12)]{IIN16FDW}). The value of a single~$e_{\alpha_1}^1(\alpha_2)$ does not change if one replaces $\Psi_{\rho_j}$ by $\Psi_{\rho_j}^\tau$, $\tau=\tau(\rho_j,x_{\alpha_2}^2)$, and $\delta_1^j$ by $\delta_1^j(\alpha_2)$, chosen appropriately. (Since $\sup|\Psi_{\rho_j}^\tau|$ is independent of~$\tau$, the diameter $\delta_1^j(\alpha_2)$ depends only on $\sup|\nabla^\ast\Psi_{\rho_j}^\tau|$. It is small for wavelets with a big angular resolution, and it is big for wavelets with a small angular resolution.) Consequently, the discrete set of wavelets constructed according to the above description is a frame for the whole~$\mathcal L^2(\mathcal S^2)$ (the frame condition is satisfied independently of~$f$), but its choice is adapted to a concrete function.

\end{document}